
\magnification1200
\font\tenBbb=msbm10 
\textfont8=\tenBbb 
\def\Bbb{\fam=8} 

\def\R{{\Bbb R}}

\def\Z{{\Bbb Z}}

\font\title=cmbx12 scaled \magstep 2
\font\bfbig=cmbx12

\def\det{\mathop{\rm det}}

\headline{\ifnum\pageno>1
   {\rlap{\sevenrm Gaetano Zampieri}
  \hfill\sevenrm Nonholonomic versus vakonomic dynamics}\else\fi\hfill
  \llap{\sevenrm }}

\centerline{\title 
Nonholonomic versus vakonomic dynamics}

\bigskip\bigskip

\centerline{Gaetano ZAMPIERI\footnote{*}{
Supported
by the MURST and by the GNFM of the CNR.}}
\bigskip

\centerline{Dipartimento di Matematica,
Universit\`a di Torino}

\centerline{via Carlo Alberto 10, 
I-10123 TORINO, Italy}

\centerline{e-mail: {\tt zampieri@dm.unito.it}}

\bigskip\bigskip

\vfill

{\narrower\smallskip\noindent
{\bf Abstract.} The main aim of the present
paper is to raise the doubt that
vakonomic dynamics  may not be
satisfactory as a model for velocity
dependent constraints. 
  
\smallskip
}
\bigskip
{\narrower\smallskip\noindent
Nonholonomic dynamics.
Variational axiomatic kind dynamics.
\smallskip
}
\bigskip

\vfill\vfill

\eject


\noindent{\bfbig 1. Introduction.} 

\bigskip

Mechanical systems with constraints on the velocities, called ``nonholonomic
constraints'', have paramount importance in engineering, in particular in
robotics, vehicular dynamics and motion generation, so they are actively and
deeply studied also from the theoretical point of view, see the recent
Bloch-Krishnaprasad-Marsden-Murray [BK], the references therein, and
Kupka-Oliva~[KO].  We are
not going to review the story of these studies, which started long ago
 with the
work of Lagrange, let us  just mention Chetaev among the most important
contributors.

The mathematical model kept essentially unchanged until about 20
years ago when a new  dynamics of 
 velocity constrained mechanical systems  was introduced 
by Kozlov~[K], and was reported in the beautiful
``Encyclopaedia of Mathematical Sciences'' 
edited by Arnol'd, see [A] Ch.~1, \S 4.
This new mechanics was called ``vakonomic'' being 
``variational axiomatic kind''.

The paper Lewis-Murray~[LM] studies the two mechanics  from the
theoretical point of view, and  deals with a ball on a rotating table
analytically, numerically, and experimentally. The experiment supports
the nonholonomic framework, moreover the authors say at page 808:
``we were not able to produce any vakonomic simulations which resembled
the experimental observations...''. However, they also say at page 809:
``Certainly, a more careful and exhaustive experimental effort on systems
other than a ball on the rotating table would be valuable in providing data
which would allow for a fair comparison on the nonholomic and vakonomic
methods.''

The previous paper Kozlov~[K], part III \S 5, relates the two
different theories to the different ways the constraints can be realized,
e.g. by large viscosity or additional masses, and it
says ``...vakonomic dynamics, which is an internally consistent model that
can be applied to the description of the motion of any mechanical systems,
is as `true' as traditional nonholonomic mechanics. The issue of the choice
of model for each {\it particular} case is ultimately resolved by
experiment.''

The aim of the present paper is to suggest the contrary: 
perhaps vakonomic mechanics is not satisfactory as a model
for  velocity dependent constraints. 
This opinion is based on the main
example used by Kozlov to support his dynamics, namely
the  skate on an inclined plane, whose nonholonomic 
behaviour is rated  ``paradoxical'' in [A] p.~19, and, at page 36,
to be compared with  the vakonomic motion he studies.

The paradox seems to come from the fact that for the nonholonomic
dynamics
``...on the average the skate
does not slide down the inclined plane...''. 

In Section~2 we recall the dynamics of {\it natural systems} with
nonholonomic  constraints, that is {\it Lagrange equation
with multipliers}, and we put it in normal form. It is a conservative
dynamics. Moreover, we see that it is {\it reversible}, 
namely the set of
solutions is invariant under time inversion. 
The example of the nonholonomic
skate ends the section. 

Section~3  deals with the same systems 
but with vakonomic dynamics. Now we have many more
solutions since new ``latent variables'' $\lambda$ appear 
([A]~p~37, and the paper [K] III at page 44, loosely
relate this fact with the ``unobservable
 quantities...in quantum mechanics'').
This situation, quite strange for classical mechanics,
is not overcome by  passing
to the Hamiltonian framework. Anyway, we prefer to keep  our
discussion in the Lagrangian framework, which is, perhaps, clearer for our
purposes.

We show that also the vakonomic dynamics of natural systems
is reversible, at least in a suitable sense, and we see no
reasons why it should not be reversible, 
as  nonholonomic dynamics is. 
The section ends with  the study of some  vakonomic solutions to
the skate  which seem paradoxical to me {\it for all initial values of the
``latent variable''}.

Kharlamov [Kh] considers also the skate on a plane 
to criticize the vakonomic
mechanics. His plane is horizontal, so he is able to obtain explicit
expressions for the vakonomic solutions which correspond to the nonholomic
uniform motions along straight lines. The author concludes that the 
vakonomic skate has a 
``fanciful track''.  However, the  uniform motions along straight lines
are  vakonomic motions. This fact does not ``save the situation'' according
to Kharlamov  since ``there is no basis on which to assign a
specific value'' to the latent variable at the initial time.

Finally, Section 4 deals with more general nonholonomic systems
and revisits the nonholonomic skate by adding some dissipation
on the rotational degree of freedom only. Now, the skate slides
eventually down. 

\goodbreak
\bigskip\bigskip\bigbreak

\noindent {\bfbig 2. Nonholonomic dynamics of natural systems.} 
 \bigskip

The  {\it Lagrange equation} of Classical Mechanics is the celebrated
$${d\over {dt}}\;{\partial L\over {\partial \dot q}}(t,q,\dot q)\,-\,
   {\partial L\over{\partial q}}(t,q,\dot q)\,=\,0\,.
\eqno(2.1)$$
The  Lagrangian function $L$ is assumed of class $C^2$ 
on some open connected
subset of $\R\times\R^N\times\R^N$.

In particular we are interested in the 
{\it natural Lagrangian functions}:
$$L(t,q,\dot q)={1\over 2} \dot q\cdot A(q) \dot
q-U(q)\,,\qquad A(q)=A(q)^T>0\,,\eqno(2.2)$$
where $A(q)^T$ is the transpose matrix of the $N\times N$ positive
definite matrix
$A(q)$, the central dot is the usual scalar product, ${1\over 2}
\dot q\cdot A(q) \dot q$ is called the kinetic energy,
 and $U(q)$ the
potential energy.

If one asks the mechanical system to obey the  {\it nonholonomic
constraint equation} 
$$B(q)\,\dot q=0\eqno(2.3)$$
where $B(q)$ is an $n\times N$, with $n<N$, {\it full rank} matrix
 at each $q$,
then one has a  nonholonomic system, which we briefly call
{\it ``nonholonomic natural system''}, whose dynamics
is ruled by the constraint (2.3) together with
the following classic {\it Lagrange equation with multipliers} 
$${d\over {dt}}\;{\partial L\over {\partial \dot q}}(t,q,\dot q)\,-\,
   {\partial L\over{\partial q}}(t,q,\dot q)\,=\,B(q)^T \mu\,.
\eqno(2.4)$$
The multipliers are the components of the new unknown function
$\mu(t)$ with values in $\R^n$. In the sequel, we assume $B\in C^2$ as
$L$. For natural systems the dynamics is then ruled by equations
of the following form
$$\cases{A(q)\,\big(\ddot q\,+
\,\Gamma(q)[\dot q, \dot q]\big)\,+\nabla U(q)\,
=\,B(q)^T \mu\cr
\noalign{\smallskip}\cr
B(q)\,\dot q\,=\,0\cr}
\eqno(2.5)$$
where $\nabla$ is the gradient operator, and
$\Gamma(q):\R^N\times\R^N\to\R^N, (u,v)\mapsto
\Gamma(q)[u,v]$, is a bilinear symmetric map whose components are  called
Christoffel symbols in the classical books. 

We easily check that it is a conservative dynamics, namely the
{\it total energy}  $\ {1\over 2} \dot q\cdot A(q) \dot
q+U(q)$ is a first integral of (2.5).

For brevity, let us omit the functional dependences for a while. The
$N\times N$ matrix $A$ is symmetric and  positive definite, in particular
then invertible, so we can define the new $n\times n$ matrix
$C:=B\,A^{-1}\,B^T$ which is nonsingular, indeed 
$\ C\,\xi=0\ \Longrightarrow\  0=\xi\cdot C\,\xi=\xi\cdot B\,A^{-1}\,B^T
\xi=B^T \xi\cdot A^{-1}\,B^T \xi\ \Longrightarrow\  B^T \xi=0\ 
$
since $A^{-1}$ is positive definite as $A$; finally $B^T \xi=0
\ \Longrightarrow\  \xi=0$ since $B$ has full rank. 

It is easy to see that the following $\,(N+n)\times (N+n)\,$ matrices 
$$\left(\matrix{A & -B^T\cr
                -B & 0\cr}\right)\,,\qquad\qquad 
\left(\matrix{A^{-1}-A^{-1} B^T C^{-1} B A^{-1} & -A^{-1} B^T C^{-1}\cr
                -C^{-1} B A^{-1} & -C^{-1}\cr}\right) \eqno(2.6)
$$
are inverse of each other.
Thus the following system in the unknown functions $q(t),\mu(t)$, which is
obtained from  (2.5) by  differentiating 
the second equation with respect to $t$,
$$\left(\matrix{A & -B^T\cr
                -B & 0\cr}\right) \left(\matrix{\ddot q\cr
                \mu\cr}\right)+\left(\matrix{A(q)\Gamma(q)[\dot q, \dot
q]+\nabla U(q)\cr
                -B'(q)[\dot q,\dot q]\cr}\right)
\,=\,0\eqno(2.7)$$
 is equivalent to
$$\left(\matrix{\ddot q\cr
                \mu\cr}\right)=
\left(\matrix{A^{-1}-A^{-1} B^T C^{-1} B A^{-1} & -A^{-1} B^T C^{-1}\cr
                -C^{-1} B A^{-1} & -C^{-1}\cr}\right)
\left(\matrix{-A(q)\Gamma(q)[\dot q, \dot q]-\nabla U(q)\cr
                B'(q)[\dot q,\dot q]\cr}\right).\eqno(2.8)$$
Reinstating the
functional dependences and  defining suitable new functions $f,g$ we get
the following ``normal form'' which is quadratic in
$\dot q$
$$\cases{\ddot q=f(q)+g(q)[\dot q,\dot q]\cr
\noalign{\smallskip}\cr
\mu=h(q)+k(q)[\dot q,\dot q]\cr}\eqno(2.9)$$
We can get rid of the last equation and the unknown $\mu$.

Since the second 
equation in (2.7) is obtained from  (2.3) by differentiation, we have 
$\,B(q)\,\dot q=const$ along the solutions of the first equation in (2.9).
So, the constraint equation $B(q)\,\dot q=0$ is satisfied provided it holds
at  some time $t_0$. Therefore, we arrive at the conclusion that {\it all
nonholonomic motions can be obtained from  Cauchy problems of the following
kind}  where an additional
 condition  restricts the choice of the initial velocity

$$\cases{\ddot q=f(q)+g(q)[\dot q,\dot q]\cr
\noalign{\smallskip}\cr
\big(q(t_0)\,,\dot q(t_0)\big)=(q_0,\dot q_0)\in \{(u,v): B(u)\,v=0\}\cr
}\eqno(2.10)$$

If $q(t)$ is a solution to (2.10),  then the  function $r(t)=q(-t)$
is also a nonholonomic motion which has initial data $(r(-t_0),\dot
r(-t_0))=(q_0,-\dot q_0)$ at time $-t_0$. So {\it the set of nonholonomic
motions is invariant under time reversal, and we can say that
 the nonholonomic
dynamics of natural systems is reversible}.

\bigskip

{\bf Example: the nonholonomic skate}. Consider a skate
 (an homogeneous material segment) on
an inclined plane with Cartesian coordinates
 $x,y$. The $x$-axis points downward
while the
$y$-axis is horizontal. $x$ and $y$ will denote the 
coordinates of the center of the skate
which has another Lagrangian coordinate: 
the rotation angle $\phi$ it makes
with the unit vector of the
$x$-axis. We have a  natural Lagrangian function  as in formula
(2.2) 
$$L(t,x,y,\phi,\dot x,\dot y,\dot \phi)=
{1\over 2}\left(\dot x^2+\dot y^2+\dot
\phi^2\right)+x\,.\eqno(2.11)
$$

Suppose the center of the skate can only have velocities parallel to it, 
namely
we consider the constraint equation
$$\dot x\, \sin \phi\, -\, \dot y\, \cos \phi\,=\,0\,.\eqno(2.12)$$

The nonholonomic  equations (2.5) in the actual case are
$$\cases{\ddot x-1= \mu\, \sin\phi\cr
         \ddot y=- \mu\, \cos\phi\cr
         \ddot \phi=0\cr
         \dot x\, \sin \phi\, -\, \dot y\, \cos \phi=0\cr}\eqno(2.13)
$$
We easily get the following Cauchy problem which particularizes (2.10)
to our system, to the initial time $t_0=0$, and to some special
initial conditions (considered in [A] p.19)
$$\cases{\ddot x=\cos^2\phi-\dot\phi\,\sin\phi\left(\dot x\,
\cos\phi+\dot y\, \sin\phi\right)\cr
         \ddot y=\cos\phi \sin\phi+\dot\phi\,\cos\phi\left(\dot x\,
\cos\phi+\dot y\,
\sin\phi\right)\cr
         \ddot \phi=0\cr
       (x(0),y(0),\phi(0),\dot x(0),\dot y(0),\dot
\phi(0))=(0,0,0,0,0,w)\cr
}\eqno(2.14)
$$
with $w\ne 0$. Notice that the initial data above are acceptable since they
satisfy the constraint at $t=0$:  $\dot x(0)
\,\sin
\phi(0)-\dot y(0)\,\cos
\phi(0)=0$. The solution is easily found
$$x(t)={1\over{2 w^2}}\,\sin^2 w t\,,\quad
  y(t)={1\over{2 w^2}}\,\left(w t-{1\over 2}\sin 2 w t\right)\,,\quad
\phi(t)= w t\,.\eqno(2.15)$$
It is a cycloid.
As remarked in [A] p.19, on the average the skate does not slide
down the inclined plane: $0\le x(t)\le 1/2 w^2$.

We are also interested in the motion starting with the conditions 
$\, (0,0,\pi/2,0,\dot y_0,0)$, which also satisfy the constraint 
$$x(t)=0\,,\qquad y(t)=\dot y_0\,t\,,\qquad \phi(t)=\pi/2\,.\eqno(2.16)$$

\goodbreak
\bigskip\bigskip\bigbreak


\noindent {\bfbig 3. Vakonomic dynamics of natural systems.} 
 \bigskip

As in Section~2, we consider the natural Lagrangian  $\; L(t,q,\dot
q)$ in (2.2), and the  constraint (2.3) with a full rank $n\times N$ matrix
$B(q)$. As is well known, the nonholonomic system (2.5) is a consequence
of d'Alembert's principle which is not  variational.
If we adopt a variational approach by requiring the motion to be 
 a stationary curve of the action functional $\;
q(\cdot)\mapsto\int_{t_0}^{t_1}L(t,q(t),\dot q(t)) dt$ among all
curves
having the same end points and satisfying the nonholomic constraints,
then we get a vakonomic motion (see [A] pp. 32, 33, 34 and [KO]
for details). Moreover, the motion $t\mapsto q(t)$ is vakonomic if and
only if  there exists a smooth curve $t\mapsto \lambda(t)$ in
$\R^n$, defined on the same time interval,  such that the pair
$(q(\cdot),\lambda(\cdot))$ is a solution to  the (unconstrained) variational
problem associated to the Lagrangian function
$$\L(t,q,\dot q,\lambda,\dot \lambda)=L(t,q,\dot q)\,-\,\lambda\cdot
B(q)\,\dot q\eqno(3.1)$$
(the scalar product being in $\R^n$ as the values of the new unknown 
$\lambda(t)$). Namely, the vakonomic motions can be obtained by the
Lagrange equations
$$\cases{\displaystyle{{d\over {dt}}\;{\partial \L\over {\partial \dot
q}}\,-\,
   {\partial \L\over{\partial q}}\,=\,0}\cr
        \displaystyle{{d\over {dt}}\;{\partial \L\over {\partial \dot
\lambda}}\,-\,
   {\partial \L\over{\partial \lambda}}\,=\,0}\cr}
\eqno(3.2)$$
As in the previous section we speak of {\it natural systems}, now {\it
vakonomic}. Their dynamics is ruled by 
$$\cases{A(q)\,\big(\ddot q\,+\,\Gamma(q)[\dot q, \dot q]\big)\,+\nabla
U(q)\,-\,B(q)^T\,\dot
\lambda\,-\lambda\cdot B'(q)\,\dot q\,+\,{\partial\over{\partial
q}}\left(\lambda\cdot B(q)\dot q\right)\,=\,0\cr
\noalign{\smallskip}\cr
B(q)\,\dot q=0\cr}
\eqno(3.3)$$
where, for each  $i\in\{1,\dots,N\}$, the $i$-components of two of the
previous expressions are
$$\left(\lambda\cdot B'(q)\,\dot
q\right)_i=\sum_{\alpha=1}^n\,\lambda_{\alpha}\,\sum_{j=1}^N\,
{\partial B_{\alpha i}\over{\partial q_j}}\;\dot q_j\,,\quad
\left({\partial\over{\partial
q}}\left(\lambda\cdot B(q)\dot q\right)\right)_i=\sum_{\alpha=1}^n\,
\lambda_{\alpha}\,\sum_{j=1}^N\,
{\partial B_{\alpha j}\over{\partial q_i}}\;\dot q_j\,.\eqno(3.4)
$$

As in Section~2 we can derive the second equation in (3.3) 
with respect to $t$
and write the system obtained from (3.3) in this way as
$$\left(\matrix{A(q) & -B(q)^T\cr
                 -B(q) & 0\cr}\right)\left(\matrix{\ddot q\cr
\dot\lambda\cr}\right)\,+\,
F(q,\dot q,\lambda)\,=\,0\eqno(3.5)$$
for a suitable $F$. Again the $(N+n)\times (N+n)$ matrix is invertible since
$A(q)=A(q)^T>0$ and $B(q)$ has full rank at
any point $q$, and the inverse is as in (2.6). So we can solve (3.5)
with respect to $(\ddot q,\dot \lambda)$ and get an equivalent system
in normal form
$$\left(\matrix{\ddot q\cr
\dot\lambda\cr}\right)\,=\,
G(q,\dot q,\lambda).\eqno(3.6)$$
 The function
$B(q) \dot q$ is a first integral, i.e. constant along the solutions. If  
$B(q(t_0))\dot q(t_0)=0$ at the initial time
$t_0$, then we have a vakonomic motion,
 moreover all vakonomic motions can be
obtained in this way. However, now $\lambda(t_0)$ is arbitrary in
the Cauchy problem. It seems that we have too many vakonomic motions. 

As in the nonholonomic dynamics of natural systems the total energy
is conserved. However, now we cannot get rid of $\lambda$. 

In order to discuss the time reversibility, let us consider the new unknown
function $\sigma(t)$, with
$\dot \sigma(t)=\lambda(t)$. Then Kozlov function (3.1) becomes quadratic
in $(\dot q,\dot \sigma)$
$$\L(t,q,\dot q,\sigma,\dot \sigma)={1\over 2} \dot q\cdot A(q) \dot
q-U(q)-\dot \sigma\cdot B(q)\dot
q\eqno(3.7)$$
and the related equations
$$\cases{\displaystyle{{d\over {dt}}\;{\partial \L\over {\partial \dot
q}}\,-\,
   {\partial \L\over{\partial q}}\,=\,0}\cr
         \displaystyle{{d\over {dt}}\;{\partial \L\over {\partial \dot
\sigma}}\,-\,
   {\partial \L\over{\partial \sigma}}\,=\,0}\cr}
\eqno(3.8)$$
give directly (3.5) with $\ddot \sigma$ instead of $\dot \lambda$, $\dot
\sigma$ instead of $\lambda$, and no dependence on $\sigma$
$$\left(\matrix{A(q) & -B(q)^T\cr
                 -B(q) & 0\cr}\right)\left(\matrix{\ddot q\cr
\ddot\sigma\cr}\right)\,+\,
F(q,\dot q,\dot\sigma)\,=\,0\eqno(3.9)$$
This system can again be put in normal form
$$\left(\matrix{\ddot q\cr
\ddot\sigma\cr}\right)\,=\,
G(q,\dot q,\dot\sigma)\eqno(3.10)$$
with the first integral $B(q)\,\dot q$. Whenever  this first integral
vanishes, we  get a vakonomic motion and all vakonomic motions can be 
obtained in this way. Finally, the solutions to (3.3) are the pairs
$\,(q(t),\lambda(t))=(q(t),\dot \sigma(t))$.

The equation (3.10) is quadratic in the velocities, thus for any solution
we have a corresponding reversed one: if $(q(t),\lambda(t))$ is a 
vakonomic motion, then so is $(q(-t),-\lambda(-t))$. 
Perhaps one  can   argue that these  solutions are not necessarily {\it
physically} acceptable at the same time, so accepting one and rejecting the
other. I think that this point should be somehow justified: without
a good reason  why should the  vakonomic dynamics of
natural systems miss the  feature of time reversibility? 
Would not the very lack of time reversibility  be a sufficient
reason to reject the vakonomic model?

\bigskip

{\bf Example: the vakonomic skate}. Consider again the skate as in
Section~2, now vakonomic. The above Kozlov function (3.1) for this
mechanical system is
$$\L(t,x,y,\phi,\dot x,\dot y,\dot \phi,\lambda,\dot\lambda)=
{1\over 2}\left(\dot
x^2+\dot y^2+\dot
\phi^2\right)+x-\lambda\,\left(\dot x\, \sin \phi\, -\, \dot y\, \cos
\phi\right)\,.\eqno(3.11)
$$
The vakonomic equations are
$$\cases{\ddot x-\dot \lambda\, \sin\phi-\lambda\, \dot\phi\, \cos\phi=1\cr
         \ddot y+\dot \lambda\, \cos\phi-\lambda \, \dot\phi\, \sin\phi=0\cr
         \ddot \phi+\lambda\left(\dot x\, \cos\phi+\dot y\,
\sin\phi\right)=0\cr
         \dot x\, \sin \phi\, -\, \dot y\, \cos \phi=0\cr}\eqno(3.12)
$$
We easily get the following system
$$\cases{\ddot x=-\dot\phi\,\sin\phi\left(\dot x\, \cos\phi+\dot y\,
\sin\phi\right)+\lambda\,\dot\phi\,\cos\phi+\cos^2\phi\cr
         \ddot y=+\dot\phi\,\cos\phi\left(\dot x\, \cos\phi+\dot y\,
\sin\phi\right)+\lambda\,\dot\phi\,\sin\phi+\sin\phi\,\cos\phi\cr
         \ddot \phi=-\lambda\left(\dot x\, \cos\phi+\dot y\,
\sin\phi\right)\cr
        \dot\lambda=-\sin\phi-\dot\phi\left(\dot x\, \cos\phi+\dot y\,
\sin\phi\right)\cr}\eqno(3.13)
$$
which has the first integral
$$  \dot x\, \sin \phi\, -\, \dot y\, \cos \phi\,=a.\eqno(3.14)$$
The solutions to (3.13) with $a=0$ are all the
vakonomic motions.

The vakonomic motions to be compared to the nonholonomic motions (2.15)
(we use the plural since there is a parameter $w$),
according to Kozlov have $\lambda(0)=0$ (this choice is probably justified
by the results it gives), and is studied in [K] III, Section~3. Their 
features are also reported in [A] p.~35, 36. Here, let us just
quote from there that the skate slides monotonically down for $t>0$,
and almost all solutions tend to turn sideways: $\phi(t)$ 
converges to one of
the points $ \pi/2+m \pi$ ($m\in\Z$), as $t\to+\infty$.
Let us add that for $t<0$ it slides monotonically up.

Now, let us go to the nonholonomic motions (2.16). First of
all, let us immediately remark that they have no corresponding
vakonomic motions the extra variable $\lambda$ notwithstanding, unless
$\dot y_0=0$,  the ``vakonomic equilibrium'',
for which the ``latent variable'' $\lambda$ has the precise
meaning of reversed time:
$$x(t)=0\,, \qquad y(t)=0\,,\qquad \phi(t)=\pi/2\,,\qquad
\lambda(t)=\lambda(0)-t\,.\eqno(3.15)$$

Consider now  (3.13) with the initial conditions $$
       (x(0),y(0),\phi(0),\dot x(0),\dot y(0),\dot
\phi(0),\lambda(0)) = (0,0,\pi/2,0,\dot y_0,0,\lambda_0)\,,\qquad \dot y_0\ne
0\,.
\eqno(3.16)$$ For the moment assume that $\,\lambda_0\ne 0$ too.
Then (3.13) give $\,\ddot x(0)=0$, the second and the third derivative of
the  functions $\phi(t)$ and $x(t)$ respectively, at $t=0$, are 
 $$\ddot \phi(0)=-\lambda_0\,\dot y_0\,\sin(\pi/2)=-\lambda_0\,\dot
y_0\,,\qquad x^{(3)}(0)=-\ddot
\phi(0)
\dot y(0) \sin^2(\pi/2)=\lambda_0\,\dot y_0^2\,.\eqno(3.17)$$
Therefore
$$\lambda_0\ne 0\qquad \Longrightarrow\qquad \cases{x(t)=\lambda_0 
\,\dot y_0^2
\, t^3\,/3!
+o(t^3)\cr
\noalign{\smallskip}\cr
y(t)=\dot y_0\,t+o(t)\cr
\noalign{\smallskip}\cr
\phi(t)=\pi/2-\lambda_0\,\dot
y_0\,t^2\,/2+o(t^2)\cr
}\eqno(3.18)$$ If we accept the reversibility 
of vakonomic dynamics we arrive
at the conclusion that either the solution we are analyzing, or
the one obtained by time reversal, {\it slides  monotonically up in
a full neighbourhood of $t=0$}. I guess we can rate this behaviour
paradoxical given the initial conditions above.

Next,
consider the other possibility $\lambda_0=0$. In this case
(3.13) gives
 $$\eqalign{\ddot \phi(0)=0,\quad & \phi^{(3)}(0)=-\dot
\lambda(0)\,\dot y(0)\,\sin(\pi/2)=\dot y_0 \,\sin^2(\pi/2)=\dot y_0, \cr
x^{(2)}(0)=0, \quad & x^{(3)}(0)=0, 
\quad x^{(4)}(0)=-\phi^{(3)}(0)\,\dot y(0)
\sin^2(\pi/2)=-\dot y_0^2\,.\cr}\eqno(3.19)$$
So $$\lambda_0= 0\qquad \Longrightarrow\qquad \cases{
x(t)=-\dot
y_0^2\,t^4\,/4!+o(t^4)\cr
\noalign{\smallskip}\cr
y(t)=\dot y_0\,t+o(t)\cr
\noalign{\smallskip}\cr
\phi(t)=\pi/2+\dot
y_0\,t^3\,/3!+o(t^3)\cr
 }\eqno(3.20)$$
which seems another paradox:  {\it the skate
slides  up
for $t>0$} (small enough) {\it while both $\;\dot x(0)$ and $\dot
\phi(0)\,$ vanish, and
$\,\lambda(0)=0\,$ too.}

Of course, also the ``body-plus-fluid'', considered in [K2] p. 598 to
answer the criticism in [Kh],  has the same strange behaviour as well
as any other physical system we may think should obey the same dynamics.

\goodbreak
\bigskip\bigskip\bigbreak


\noindent {\bfbig 4. General nonholonomic dynamics.} 
 \bigskip

For the final example, equations only slightly more general than
those considered in Section~2 are needed. But full generality can
be obtained with only little work from what described in Section~2,
so we proceed to it.

The Lagrange
equation whith multipliers (2.4) is also  studied for general Lagrangian
functions
$L(t,q,\dot q)$, and  has been also extended to constraints possibly
non-linear in
$\dot q$ by the work of Chetaev and others. 
Notice that  justifying the extension was
a difficult job, see the recent   Cardin-Favretti~[CF],
and  Massa-Pagani~[MP];  but concrete
non-linear examples can be given,  see  Benenti~[B].

The general equations for {\it nonholonomic dynamics} are 
$$\cases{\displaystyle{{d\over {dt}}\;{\partial L\over {\partial \dot
q}}(t,q,\dot q)\,-\,
   {\partial L\over{\partial q}}(t,q,\dot q)\,=\,\left({\partial
b\over{\partial \dot q}}(t,q,\dot q)\right)^T\,\mu}
\cr 
\noalign{\smallskip}\cr
b(t,q,\dot q)=0\cr}
\eqno(4.1)$$
where, for each $t$, the values $q(t),\dot q(t)\in\R^N$, $\mu(t),
b(t,q(t),\dot q(t))\in\R^n$. The standard conditions to put (4.1) in normal
form are
$$\det \,S(t,q,\dot
q)\ne 0\,\qquad \hbox{\rm where}\quad \,S(t,q,\dot
q):={\partial^2 L\over {\partial
\dot q^2}}(t,q,\dot q)\,\qquad \hbox{\rm at all}\quad (t,q,\dot
q)\,,\eqno(4.2)$$
$$\eqalign{ &\hbox{\rm and}\quad \det R(t,q,\dot q)\ne 0\qquad \hbox{\rm at
all}\quad (t,q,\dot q)\,,\cr
 &\hbox{\rm where}\quad 
R:=D\,S^{-1}\,D^T\,,\quad  \hbox{\rm with}\quad D(t,q,\dot q):={\partial
b\over{\partial
\dot q}}(t,q,\dot q)\,.\cr}\eqno(4.3)$$

Indeed, differentiating the second equation (4.1) with respect to $t$,
 and introducing suitable new functions  $r$ and $s$, the system
(4.1) gives
$$\cases{S(t,q,\dot q)\,\ddot q-r(t,q,\dot q)\,=\,D(t,q,\dot q)^T\,\mu 
\cr 
\noalign{\smallskip}\cr
-D(t,q,\dot q)\,\ddot q-s(t,q,\dot q)\,=\,0\cr}
\eqno(4.4)$$

Next, we check at once that  the following two $\,(N+n)\times (N+n)\,$
matrices are inverse of each other 
$$\left(\matrix{S & -D^T\cr
                -D & 0\cr}\right)\,,\qquad\qquad 
\left(\matrix{S^{-1}-S^{-1} D^T R^{-1} D S^{-1} & -S^{-1} D^T R^{-1}\cr
                -R^{-1} D S^{-1} & -R^{-1}\cr}\right)\,. \eqno(4.5)
$$
Thus the system (4.4) is equivalent to the following one in normal form
$$\left(\matrix{\ddot q\cr
                \mu\cr}\right)\,=\,
\left(\matrix{S^{-1}-S^{-1} D^T R^{-1} D S^{-1} & -S^{-1} D^T R^{-1}\cr
                -R^{-1} D S^{-1} & -R^{-1}\cr}\right)\, \left(\matrix{r\cr
                s\cr}\right)\,.\eqno(4.6)$$

 Since the second 
equation in (4.4) is obtained from the constraint equation by derivation
in $t$, we have that $b(t,q,\dot q)=const$ along the solutions of the
first equation in (4.6).
So,  $b(t,q,\dot q)=0$ is satisfied provided
it holds at  the initial time $t_0$. Finally, we can get rid of the
multipliers  arriving at the following  Cauchy problem which generalizes
(2.10) 
$$\cases{\ddot q=F(t,q,\dot q)\cr
\noalign{\smallskip}\cr
\big(q(t_0)\,,\dot q(t_0)\big)=(q_0,\dot q_0)\in \{(u,v): b(t_0,u,v)=0\}\cr
}\eqno(4.7)$$
where the first equation
is the first equation in (4.6),  briefly written. 

\bigskip

{\bf Example: the nonholonomic skate with rotational friction}.
We modify the Lagrangian function (2.11)  so it has rotational
dissipation 
$$L(t,x,y,\phi,\dot x,\dot y,\dot \phi)={1\over 2}\left(\dot x^2+\dot
y^2+e^{k t}\,\dot
\phi^2\right)+x\eqno(4.8)
$$
where $k> 0$ is a new parameter.

The nonholonomic  equations are
$$\cases{\ddot x-1= \mu\, \sin\phi\cr
         \ddot y=- \mu\, \cos\phi\cr
         \ddot \phi+k \dot \phi=0\cr
         \dot x\, \sin \phi\, -\, \dot y\, \cos \phi=0\cr}\eqno(4.9)
$$
Only the third equation is different from the one in (2.13). Instead of
(2.14) we then have
$$\cases{\ddot x=\cos^2\phi-\dot\phi\,\sin\phi\left(\dot x\,
\cos\phi+\dot y\, \sin\phi\right)\cr
         \ddot y=\cos\phi \sin\phi+\dot\phi\,\cos\phi\left(\dot x\,
\cos\phi+\dot y\,
\sin\phi\right)\cr
         \ddot \phi=-k \dot \phi\cr
       (x(0),y(0),\phi(0),\dot x(0),\dot y(0),\dot
\phi(0))=(0,0,0,0,0,w)\cr
}\eqno(4.10)
$$
The solution is now
$$\eqalign{\phi(t)&={w\over k}\left(1-e^{-k t}\right)\,,\qquad  x(t)={1\over
2}\left(\int_0^t\,\cos\left(\phi(\tau)\right)\,d\tau\right)^2\,,\cr
y(t)&=\int_0^t\,\sin\left(\phi(\xi)\right)\int_0^{\xi}\,
\cos\left(\phi(\tau)\right)\,d\tau d\xi\,.\cr}\eqno(4.11)$$
Notice from $x(t)$ that the skate slides eventually down. Indeed, if
$$\lim_{t\to +\infty} \phi(t)={w\over k}\notin \left\{{\pi\over 2}+m \pi:
m\in\Z\right\}\qquad
\hbox{{\rm then}}\qquad  \lim_{t\to+\infty} x(t)=+\infty\,;\eqno(4.12)$$
and otherwise $x(t)$ converges to a strictly positive finite limit
 as $t\to+\infty$.

\goodbreak
\bigskip\bigskip\bigbreak
\frenchspacing

\noindent {\bfbig  References.}

\nobreak
\bigskip
  
\item{[A]} V.I. Arnol'd (Ed.),  {\it Dynamical Systems III}. 
Encyclopaedia of Mathematical Sciences, Vol.~3. Springer-Verlag
(1988). 

\medskip\goodbreak

\item{[B]} S. Benenti, {\sl Geometrical aspects of the dynamics of non-holonomic
systems}, Rend. Sem. Mat. Univ. Pol. Torino~{\bf 54} (1996), 203--212.

\medskip\goodbreak

\item{[BK]} A.M. Bloch, P.S. Krishnaprasad, J.E. Marsden \& R.M. Murray,
{\sl Nonholonomic mechanical systems with symmetry}, Arch. Rational Mech.
Anal.~{\bf 136} (1996), 21--99.

\medskip\goodbreak
 
\item{[CF]} F. Cardin \& M. Favretti, {\sl On nonholonomic and vakonomic
dynamics of  mechanical systems with nonintegrable
constraints}, J. Geom. and Phys.~{\bf 18} (1996), 295--325.

\medskip\goodbreak

\item{[Kh]} P.V. Kharlamov, {\sl A critique of some mathematical
models of mechanical systems with differential constraints}, J. Appl.
Math. Mech. {\bf 56} (1992) 584--594.

\medskip\goodbreak

\item{[K]} V.V. Kozlov, {\sl Dynamics of systems with nonintegrable
constraints}:  Part {\sl I}, Moscow Univ. Mech. Bull. {\bf 37} (1982),
27--34; Part {\sl II},  74--80;
Part {\sl III}, vol. {\bf 38} (1983), 40--51.

\medskip\goodbreak

\item{[K2]} V.V. Kozlov, {\sl The problem of realizing constraints in
dynamics}, J. Appl. Math. Mech. {\bf 56} (1992) 594--600.

\medskip\goodbreak

\item{[KO]} I. Kupka \& W.M. Oliva, {\sl The non-holonomic mechanics},
Preprint Math. Dept. Inst. Sup. T\'ecnico, Lisbon (1997).

\medskip\goodbreak

\item{[LM]} A.D. Lewis \& R.M. Murray, {\sl Variational principles for
constrained systems: theory and experiment}, Int. J. Non-Linear
Mechanics {\bf 30} (1995), 793--815.

\medskip\goodbreak

\item{[MP]} E. Massa \& E. Pagani, {\sl Classical dynamics of non-holonomic
systems: a geometric approach}, Ann. 
Inst. H. Poincar\'e Phys. Th\'eo. {\bf 55}
(1991), 511--544.

\goodbreak

\bye